\documentclass[11pt,a4paper]{article}

\usepackage{epic,eepic,epsfig,amsmath,amssymb,amsthm}
\usepackage{t1enc}

\usepackage{color}

\usepackage{bbm}

\def\virgp{\raise 2pt\hbox{,}}

\renewcommand{\geq}{\geqslant}
\renewcommand{\leq}{\leqslant}
\def\N{{\mathbb N}}

\def\R{{\mathbb R}}

\def\virgp{\raise 2pt\hbox{,}}
\def\cdotpv{\raise 2pt\hbox{;}}

\def\1{\mathbbm{1}}

\newtheorem{theorem}{Theorem}[section]
\newtheorem{corollary}[theorem]{Corollary}
\newtheorem{proposition}[theorem]{Proposition}
\newtheorem{lemma}[theorem]{Lemma}

\theoremstyle{remark}
\newtheorem{remark}{Remark}[section]

\theoremstyle{definition}
\newtheorem{definition}{Definition}[section]

\theoremstyle{definition}

\theoremstyle{definition}

\begin{document}

\title{A simple way to get the box dimension of the graph of the~W. function }

\author{\LARGE{\textbf{Claire David}}}

\maketitle
\centerline{Sorbonne Universit\'e}

\centerline{CNRS, UMR 7598, Laboratoire Jacques-Louis Lions, 4, place Jussieu 75005, Paris, France}


\section*{Introduction}

\hskip 0.5cm
The determination of the box and Hausdorff dimension of the graph of the Weierstrass function has, since long been, a topic of interest. In the following, we show that the box-counting dimension (or Minskowski dimension) can be obtained directly, without using dynamical systems tools.\\

Let us recall that, given~$\lambda \,\in\,]0,1[$, and~$b$ such that~$\lambda\,b > 1+ \displaystyle \frac{3\,\pi}{2}$, the  Weierstrass function\index{Weierstrass}

$$ x \,\in\,\R \mapsto  { \mathcal W}( x)=\displaystyle \sum_{n=0}^{+\infty} \lambda^n\,\cos \left (   \pi\,b^n\,x \right)
$$

\noindent  is continuous everywhere, while nowhere differentiable. The original proof, by~K.~Weierstrass~\cite{Weierstrass1872}, can also be found in~ \cite{Titschmarsh1939}. It has been completed by the one, now a classical one, in the case where~$  \lambda\, b > 1$, by~G.~Hardy \cite{Hardy1911}.\\

After the works of A.~S.~Besicovitch and H.~D.~Ursell~\cite{BesicovitchUrsell1937}, it is Beno\^{i}t~Mandelbrot~\cite{Mandelbrot1977} who particularly highlighted the fractal properties of the graph of the Weierstrass function. He also conjectured that the Hausdorff dimension of the graph is~\mbox{$D_{\mathcal W}= 2 + \displaystyle \frac{\ln \lambda}{\ln b}$}. Interesting discussions in relation to this question have been given in the book of K.~Falconer \cite{Falconer1985}. A series of results for the box dimension can be found in the works of J.-L.~Kaplan et al.~\cite{Kaplan1984} (where the authors show that it is equal to the Lyapunov dimension of the equivalent attracting torus), in the one of ~F.~Przytycki and M.~Urba\'nki~\cite{PrzytyckiUrbanski1989}, and in those by T-Y.~Hu and K-S.~Lau~\cite{HuLau1993}. As for the Hausdorff dimension, a proof was given by~B.~Hunt \cite{Hunt1998} in 1998 in the case where arbitrary phases are included in each cosinusoidal term of the summation. Recently, K.~Bara\'{n}sky, B.~B\'{a}r\'{a}ny and J.~Romanowska~\cite{Baransky2014} proved that, for any value of the real number~$b$, there exists a threshold value~$\lambda_b$ belonging to the interval~\mbox{$\left ]\displaystyle \frac{1}{b},1\right [$}  such that the aforementioned dimension is equal to~\mbox{$D_{\mathcal W}$} for every~$b$ in~\mbox{$ \left ]\lambda_b,1 \right[$}. Results by W.~Shen~\cite{WShen2015} go further than the ones of~\cite{Baransky2014}. In~\cite{Keller2017}, G.~Keller proposes what appears as a much simpler and very original proof.\\

May one wish to understand the demonstrations mentioned above, it requires theoretical background in dynamic systems theory. For instance, in the work of J.-L.~Kaplan et al.~\cite{Kaplan1984}, the authors  call for results that cannot be understood without knowledge on the Lyapunov dimension. One may also note that their proof, which enables one to obtain the box-counting dimension of the aforementioned graph, involves sequences revolving around the Gamma Function, Fourier coefficients, integration in the complex plane, definition of a specific measure, the solving of several equations, thus, a lot of technical manipulations (on eleven pages), to yield  the result.  \\

Following those results,~F.~Przytycki and M.~Urba\'nki~\cite{PrzytyckiUrbanski1989} give a general method leading to the value of this box-counting dimension. It was initially devoted to the calculation of the Hausdorff dimension of the graph. It appears simpler than the one by~Kaplan et al., calling for Frostman's lemma~\cite{PrzytyckiUrbanskiBook}, ~\cite{Pesin}. The authors deal with continuous functions~$f$ satisfying conditions of the form:\\

 $ \forall\,(x_1,x_2)\,\in\, [0,1]^2 \, : $ \vskip -0.5cm
 
$$ \displaystyle \sup \, \left  \lbrace f(a_1)-f(a_2) \, \big   | \, x_1   \leq a_1 \leq  a_2\leq x_2   \right \rbrace \geq C\,  |   x_1-x_2|^\alpha \qquad (\star)   
$$

\noindent  where~$C$ and~$\alpha< 1$ denote strictly positive real constant.\\

In order to apply the results by~F.~Przytycki and M.~Urba\'nki, one thus requires the estimate~($\star$), which is not that easy to prove. The same kind of estimate is required to obtain the Hausdorff dimension of the graph. As evoked above, existing work in the literature all call for the theory of dynamical systems. \\

Until now, the simplest calculation  is the one by~G.~Keller~\cite{Keller2017}, where the author bypasses the Ledrappier-Young theory on hyperbolic measures~\cite{LedrappierYoung}, ~\cite{Ledrappier}, embedding the graph into an attractor of a dynamical system. The proof requires~$b-$baker maps, acting on the unit square.  It also requires results on stable and unstable manifolds, as well as on related fibers.\\

In our work~\cite{ClaireGB}, where we build a Laplacian on the graph of the Weierstrass function~$\mathcal W$, we came across a simpler means of computing the box dimension of the graph, using a sequence a graphs that approximate the studied one, bypassing all the aforementioned tools. The main computation, which, for any small interval~$[x_1,x_2]\subset [0,1]$, leads to an estimate of the form:

$$ C_{inf}\, |   x_1-x_2|^{2-D_{\mathcal W}}    \leq |  { \mathcal W}(x_1)-{ \mathcal W}(x_2) |  \leq  C_{sup}\, |   x_1-x_2|^{2-D_{\mathcal W}} $$

 \noindent  where $
D_{{\mathcal W}}=2+\displaystyle \frac{\ln \lambda}{\ln  b} $, and where~$C_{inf}$ and ~$C_{sup}$ denote strictly positive constants, is done in barely two pages, and does not require specific knowledge, putting the result at the disposal of a wider audience. The key results are exposed in the sequel.
\vskip 1cm

\section{Framework of the study}

  In this section, we recall results that are developed in~\cite{ClaireGB}. We consider the case when the real number~$b$ is an integer, that we this choose to denote by~$N_b$.

\vskip 1cm

\noindent \textbf{Notation 1.1}\\

\noindent We will denote by~$\N$ the set of natural integers, and set:

$$\N^\star= \N\setminus \left \lbrace 0 \right \rbrace.$$

\vskip 1cm

 \noindent {\textbf{Notation 1.2}}\\ \\
 \noindent	In the following,~$\lambda$ and~$N_b$ are two real numbers such that:
	
$$ 0 <\lambda<1 \quad, \quad N_b\,\in\,\N \quad  \text{and} \quad \lambda\,N_b > 1 .$$
	
 \noindent  We will consider the Weierstrass function~${\mathcal W}$, defined, for any real number~$x$, by:
	
$$   {\mathcal W}(x)=\displaystyle \sum_{n=0}^{+\infty} \lambda^n\,\cos \left ( 2\,  \pi\,N_b^n\,x \right).
$$
	 
\vskip 1cm

\noindent \textbf{Property 2.1.} \textbf{Periodic properties of the Weierstrass function}\\
	\noindent For any real number~$x$:
	
	$$   {\mathcal W}( x+1)=\displaystyle \sum_{n=0}^{+\infty} \lambda^n\,\cos \left ( 2\,\pi\,N_b^{n }\,x +2\,\pi\,N_b^{n }\right)
	=\displaystyle \sum_{n=0}^{+\infty} \lambda^n\,\cos \left ( 2\,\pi\,N_b^{n }\,x  \right)={\mathcal W}( x ).
	$$
	
	
	\noindent The study of the Weierstrass function can be restricted to the interval~$[0,1[$.
	
\vskip 1cm

We place ourselves, in the sequel, in the Euclidean plane of dimension~2, referred to a direct orthonormal frame. The usual Cartesian coordinates are~$(x,y)$.\\

\vskip 1cm

 The restriction~$\Gamma_{\mathcal W}$ to~$[0,1[ \times \R$, of the graph of the Weierstrass function, is approximated by means of a sequence of graphs, built through an iterative process. To this purpose, we introduce the iterated function system of the family of~$C^\infty$ contractions from~$\R^2$ to~$\R^2$:
$$ \left \lbrace T_{0},\hdots,T_{N_b-1} \right \rbrace$$

 \noindent where, for any integer~$i$ belonging to~$\left \lbrace 0,\hdots,N_b-1  \right \rbrace$, and any~$(x,y)$ of~$\R^2$:
$$ T_i(x,y) =\left( \displaystyle \frac{x+i}{N_b}, \lambda\, y + \cos\left(2\,\pi \,\left(\frac{x+i}{N_b}\right)\right) \right) $$

\vskip 1cm

 \noindent  {\textbf{Notation}}\\ \\
 \noindent  We will denote by:
	
$$ D_{{\mathcal W}}=2+\displaystyle \frac{\ln \lambda}{\ln N_b}$$
	
 \noindent  the Hausdorff dimension of~$\Gamma_{\mathcal W}$.

\vskip 1cm

\begin{proposition}
	
$$  \Gamma_{\mathcal W} =  \underset{  i=0}{\overset{N_b-1}{\bigcup}}\,T_{i}(\Gamma_{\mathcal W}). $$
	
\end{proposition}

\vskip 1cm

\begin{definition}
 For any integer~$i$ belonging to~$\left \lbrace 0,\hdots,N_b-1\right \rbrace  $, let us denote by:
	
$$P_i=(x_i,y_i)=\left(\displaystyle \frac{i}{N_b-1},\displaystyle\frac{1}{1-\lambda}\,\cos\left ( \displaystyle\frac{2\,\pi\,i}{N_b-1}\right ) \right)$$
	
 \noindent  the fixed point of the contraction~$T_i$.\\
	
  We will denote by~$V_0$ the ordered set (according to increasing abscissa), of the points:
	
$$
\left \lbrace P_{0},\hdots,P_{N_b-1}\right \rbrace .$$

  The set of points~$V_0$, where, for any~$i$ of~\mbox{$\left \lbrace  0,\hdots,N_b-2  \right \rbrace$}, the point~$P_i$ is linked to the point~$P_{i+1}$, constitutes an oriented graph (according to increasing abscissa)), that we will denote by~$ \Gamma_{{\mathcal W}_0}$.~$V_0$ is called the set of vertices of the graph~$ \Gamma_{{\mathcal W}_0}$.\\
	
  For any natural integer~$m$, we set:
$$ V_m =\underset{  i=0}{\overset{N_b-1}{\bigcup}}\, T_i \left (V_{m-1}\right ).$$

  The set of points~$V_m$, where two consecutive points are linked, is an oriented graph (according to increasing abscissa), which we will denote by~$ \Gamma_{{\mathcal W}_m}$.~$V_m$ is called the set of vertices of the graph~$ \Gamma_{{\mathcal W}_m}$. We will denote, in the sequel, by
$${\mathcal N}^{\mathcal S}_m=2\, N_b^m+  N_b-2$$
	
 \noindent  the number of vertices of the graph~$ \Gamma_{{\mathcal W}_m}$, and we will write:
	
$$ V_m = \left \lbrace P_0^m,  P_1^m, \hdots,  P_{{\mathcal N}^ {\mathcal S}_m-1}^m \right \rbrace .$$
	
\end{definition}

\begin{figure}[h]
	\centering{\psfig{height=10cm,width=12cm,angle=0,file=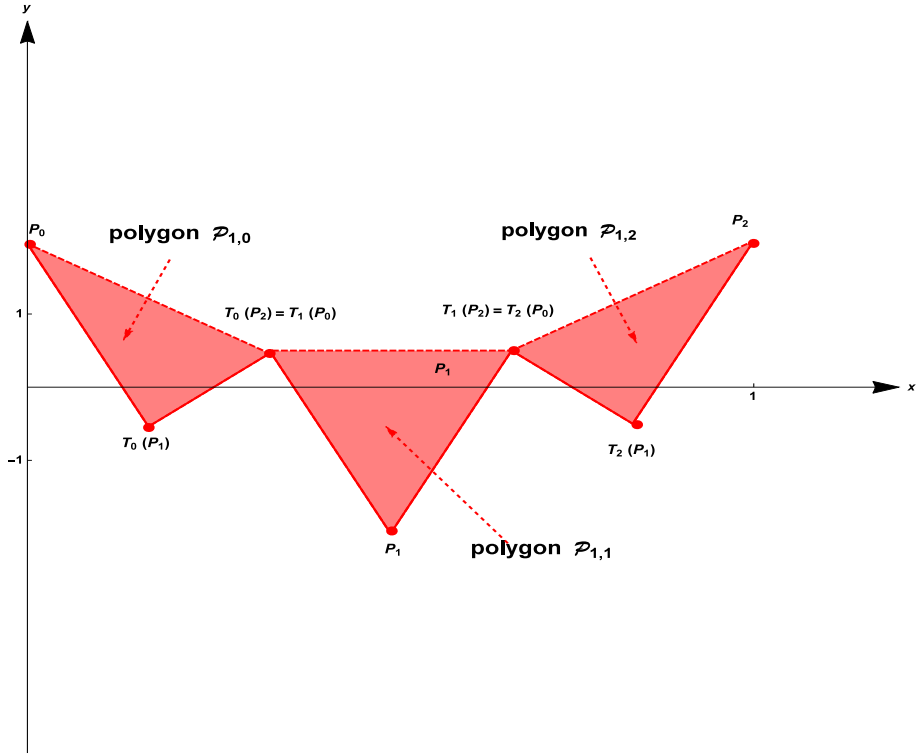} 
	\caption{The polygons~${\mathcal P}_{1,0}$,~${\mathcal P}_{1,1}$,~${\mathcal P}_{1,2}$, in the case where~$\lambda= \displaystyle \frac{1}{2}$, and~$N_b=3$.}}
\end{figure}

\begin{figure}[h]
	\centering{\psfig{height=8cm,width=12cm,angle=0,file=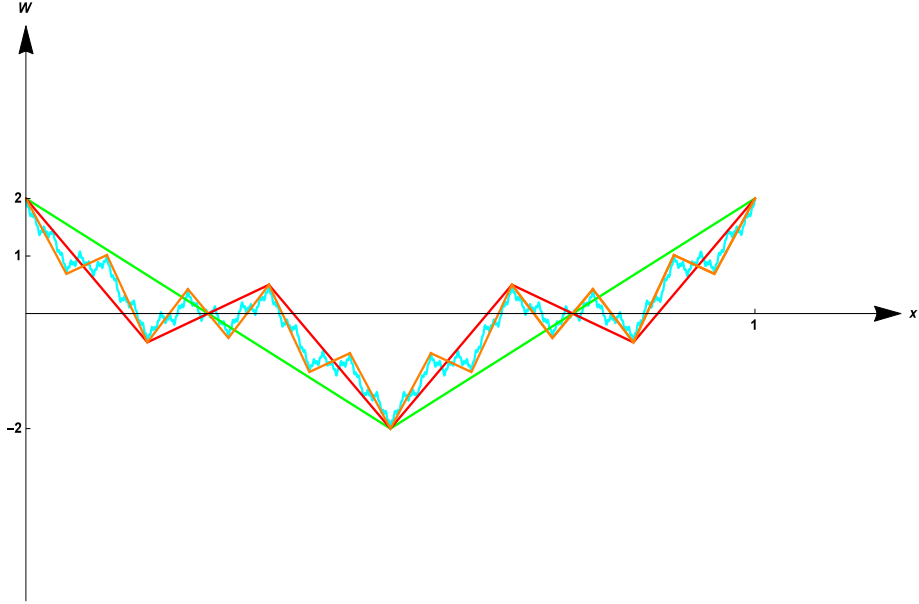} 
	\caption{The graphs~$ \Gamma_{{\mathcal W}_0 }$ (in green), ~$ \Gamma_{{\mathcal W}_1 }$ (in red),~$ \Gamma_{{\mathcal W}_2 }$ (in orange),~$ \Gamma_{{\mathcal W} }$ (in cyan), in the case where~$\lambda= \displaystyle \frac{1}{2}$, and~\mbox{$N_b=3$}.    }}
\end{figure}

\vskip 1cm

\newpage

\begin{definition}\textbf{Consecutive vertices on the graph~$\Gamma_{ \mathcal W} $ }\\
	
  Two points~$X$ and~$Y$ of~$\Gamma_{{ \mathcal W} }$ will be called \textbf{\emph{consecutive vertices}} of the graph~$\Gamma_{ \mathcal W} $ if there exists a natural integer~$m$, and an integer~$j $ of~\mbox{$\left \lbrace  0,...,N_b-2  \right \rbrace$}, such that:
	
$$ \left \lbrace \begin{array}{ccc}X &=& \left (T_{i_1}\circ \hdots \circ T_{i_m}\right)(P_j)\\
	Y &=& \left (T_{i_1}\circ \hdots \circ T_{i_m}\right)(P_{j+1})\end{array}\right.
	\qquad\left (  i_1,\hdots, i_m \right )\,\in\,\left \lbrace  0,\hdots,N_b-1  \right \rbrace^m $$
	
 \noindent	 or:

$$ \left \lbrace \begin{array}{ccc}X &=&  \left (T_{i_1}\circ  T_{i_2}\circ \hdots \circ T_{i_m}\right)\left (P_{N_b-1}\right) \\   Y &=&\left (T_{i_1+1}\circ T_{i_2}\hdots \circ T_{i_m} \right)(P_{0})\end{array}\right.  .$$

\end{definition}

\vskip 1cm

\begin{definition}
	 For any natural integer~$m$, the~$ {\mathcal N}^{\mathcal S}_m$ consecutive vertices of the graph~$  \Gamma_{{\mathcal W}_m} $ are, also, the vertices of~$N_b^m$ simple polygons~${\mathcal P}_{m,j}$,\\~\mbox{$0 \leq j \leq N_b^m-1$}, with~$N_b$ sides. For any integer~$j$ such that~\mbox{$0 \leq j \leq N_b^m-1$}, one obtains each polygon by linking the point number~$j$ to the point number~$j+1$ if~\mbox{$j = i \, \text{mod } N_b$},~\mbox{$0 \leq i \leq N_b-2$}, and the point number~$j$ to the point number~$j-N_b+1$ if~\mbox{$j =-1 \, \text{mod } N_b$}. These polygons generate a Borel set of~$\R^2$.
	
\end{definition}

\vskip 1cm
\begin{definition}\textbf{Word, on the graph~$\Gamma_{ \mathcal W} $}\\

	 Let~$m  $ be a strictly positive integer. We will call \textbf{number-letter} any integer~${\mathcal M}_i$ of~\mbox{$\left \lbrace 0, \hdots, N_b-1 \right \rbrace $}, and \textbf{word of length~$|{\mathcal M}|=m$}, on the graph~$\Gamma_{ \mathcal W} $, any set of number-letters of the form:

$$ {\mathcal M}=\left ( {\mathcal M}_1, \hdots, {\mathcal M}_m\right).$$
	
	 We will write:
	
$$ T_{\mathcal M}= T_{{\mathcal M}_1} \circ \hdots \circ  T_{{\mathcal M}_m}  .$$
	
\end{definition}

\vskip 1cm

\begin{definition}\textbf{Edge relation, on the graph~$\Gamma_{ \mathcal W} $}\\

	 Given a natural integer~$m$, two points~$X$ and~$Y$ of~$\Gamma_{{ \mathcal W}_m}$ will be called \emph{\textbf{adjacent}} if and only if~$X$ and~$Y$ are two consecutive vertices of~$\Gamma_{{ \mathcal W}_m}$. We will write:
	
$$ X \underset{m }{\sim}  Y .$$

	 This edge relation ensures the existence of a word~{${\mathcal M}=\left ( {\mathcal M}_1, \hdots, {\mathcal M}_m\right)$} of length~$ m$, such that~$X$ and~$Y$ both belong to the iterate:
	
$$T_{\mathcal M} \,V_0=\left (T_{{\mathcal M}_1} \circ \hdots \circ  T_{{\mathcal M}_m} \right) \,V_0 .$$

	 Given two points~$X$ and~$Y$ of the graph~$\Gamma_{ \mathcal W} $, we will say that~$X$ and~$Y$ are \textbf{\emph{adjacent}} if and only if there exists a natural integer~$m$ such that:
$$X  \underset{m }{\sim}  Y. $$
\end{definition}
 
\vskip 1cm
 
\begin{proposition}\textbf{Adresses, on the graph of the Weierstrass function}\\

	 Given a strictly positive integer~$m$, and a word~${\mathcal M}=\left ( {\mathcal M}_1, \hdots, {\mathcal M}_m\right)$ of length~$m\,\in\,\N^\star$, on the graph~$\Gamma_{ {\mathcal W}_m  }$, for any integer~$j$ of~$\left \lbrace 1,\hdots,N_b-2\right \rbrace  $, any~$X=T_{\mathcal M}(P_j)$ of~$ V_m \setminus V_{0}$, i.e. distinct from one of the~$N_b $ fixed point~$P_i$,
	~\mbox{$0 \leq i \leq N_b-1$}, has exactly two adjacent vertices, given by:
	
$$T_{\mathcal M}(P_{j+1})\quad \text{and} \quad T_{\mathcal M}(P_{j-1})$$
	
 \noindent	 where:
	
	$$T_{\mathcal M}  = T_{{\mathcal M}_1} \circ \hdots \circ  T_{{\mathcal M}_m}  . $$
	
	 By convention, the adjacent vertices of~$T_{{\mathcal M} }(P_{0})  $ are~$T_{{\mathcal M} }(P_{1})$ and~$T_{{\mathcal M} }(P_{N_b-1})$,
	those of~$T_{{\mathcal M} }(P_{N_b-1})$,~$T_{{\mathcal M} }(P_{N_b-2})$ and~$T_{{\mathcal M} }(P_{0})$ .

\end{proposition}

\vskip 1cm

 \noindent {\textbf{Notation 2.3.}}\\ \\
 \noindent	 For any integer~$j $ belonging to~\mbox{$\left \lbrace  0,\hdots,N_b-1  \right \rbrace $}, any natural integer~$m$, and any word~$\mathcal M$ of length~$m$, we set:
	$$\begin{array}{ccc} T_{\mathcal M}  \left (P_{j } \right) & =&
	\left ( x \left (T_{\mathcal M}  \left (P_{j } \right) \right),  y \left (T_{\mathcal M}  \left (P_{j } \right) \right) \right)
	\end{array} $$
	
	$$\begin{array}{ccc} T_{\mathcal M}  \left (P_{j+1 } \right) & =&
	\left ( x \left (T_{\mathcal M}  \left (P_{j+1 } \right) \right),  y \left (T_{\mathcal M}  \left (P_{j+1 } \right) \right) \right)
	\end{array}$$

	$$L_m= x\left ( T_{\mathcal M}  \left (P_{j+1 } \right) \right)-x\left ( T_{\mathcal M}  \left (P_{j } \right)\right)=
	\displaystyle \frac{1}{(N_b-1)\, N_b^m} $$
	
	$$h_{j,m}= y\left ( T_{\mathcal M}  \left (P_{j+1 } \right) \right)-y\left ( T_{\mathcal M}  \left (P_{j } \right)\right)  .$$

\vskip 1cm

\newpage
\section{Main results}

\begin{theorem}\textbf{An upper and lower bound, for the box-dimension of the graph~$\Gamma_{\mathcal W}$}\\
	
	 For any integer~$j$ belonging to~\mbox{$\left \lbrace 0, 1, \hdots, N_b-2 \right \rbrace $}, each natural integer~$m$, and each word~$\mathcal M$ of length~$m$, let us consider the rectangle, the width of which is:

$$L_m= x\left ( T_{\mathcal M}  \left (P_{j+1 } \right) \right)-x\left ( T_{\mathcal M}  \left (P_{j } \right)\right)=
	\displaystyle \frac{1}{(N_b-1)\, N_b^m} $$

 \noindent	 and height~$|h_{j,m}|$, such that the points~\mbox{$T_{\mathcal M}  \left (P_{j  }\right) $} and~\mbox{$T_{\mathcal M}  \left (P_{j+1 } \right )$} are two vertices of this rectangle.\\
 \noindent	 Then:\\
 
\noindent i. When the integer~$N_b$ is odd:
 	
 	\footnotesize
 	$$L_m^{2-D_{\mathcal W}}  \,(N_b-1)^{2-D_{\mathcal W}}\,  \left \lbrace   
 	\displaystyle\frac{2}{1-\lambda}\,    \sin\left ( \displaystyle\frac{ \pi }{ N_b-1 }\right )  \displaystyle \min_{0 \leq j \leq N_b-1 }  \left | \sin\left ( \displaystyle\frac{ \pi\,(2\,j+1 )}{  N_b-1 }\right )   \right| - \displaystyle  \frac{ 2\,\pi   }{  N_b \,(N_b-1)}  \,
 	\displaystyle  \frac{ 1 }{   \lambda \, N_b  -1}  \right \rbrace  \leq |h_{j,m}|  $$
 	
 	\normalsize

	 \noindent ii. When the integer~$N_b$ is even:

	 \tiny
	 $$L_m^{2-D_{\mathcal W}}   (N_b-1)^{2-D_{\mathcal W}}   \max  \left \lbrace   \! \! 
	 \displaystyle\frac{2}{1-\lambda}     \sin\left ( \displaystyle\frac{ \pi }{ N_b-1 }\right ) \! \!  
	  \displaystyle \min_{0 \leq j \leq N_b-1 } \! \!  \left | \sin\left ( \displaystyle\frac{ \pi\,(2\,j+1 )}{  N_b-1 }\right )   \right| \! \! -\! \!  \displaystyle  \frac{ 2\,\pi   }{  N_b \,(N_b-1)}   
	 \displaystyle  \frac{ 1 }{   \lambda \, N_b  -1} , 
	 \displaystyle  \frac{4   }{  N_b^{2}  }      \displaystyle  \frac{1-  N_b^{ -2}  }
	 {    N_b^2-1  }  
	  \right \rbrace  \! \! \leq \!  |h_{j,m}| \!  $$
	  
	  \normalsize
	
 \noindent	 Also:
	
	$$ |h_{j,m}| \leq  \eta_{ {\mathcal W}   }\,L_m^{2-D_{\mathcal W}} \, (N_b-1)^{2-D_{\mathcal W}}\,  $$
	
	 \noindent where the real constant~\mbox{$  \eta_{ {\mathcal W}   }$} is given by :
	
	$$ \eta_{ {\mathcal W}   }  = 2\, \pi^2\,  \left \lbrace
	\displaystyle  \frac{   (2\,N_b-1)\, \lambda\, (N_b^2-1) } {(N_b-1)^2 \, (1- \lambda )\,(\lambda \,N_b^{   2}-1) }    +
	\displaystyle  \frac{  2\, N_b } {   (\lambda \,N_b^{ 2  }-1)\, (\lambda \,N_b^{ 3  }-1)  } \right \rbrace .
$$

	

	

\end{theorem}

\vskip 1cm

\begin{corollary}
	The box-dimension of the graph~$\Gamma_{\mathcal W}$ is exactly~$D_{\mathcal W}$.

\end{corollary}

\vskip 1cm

\begin{proof}

	\noindent By definition of the box-counting dimension~$D_{\mathcal W}$~(we refer, for instance, to~\cite{Falconer1985}), the smallest number of squares, the side length of which is at most equal to~$L_m$, that can cover the graph~$\Gamma_{ \mathcal W}$ on~$[0,1[$, obeys, approximately, a power law of the form:
	
	$$c\, L_m^{-D_{\mathcal W}} \quad , \quad c >0.$$

	Let us set:\\
	\footnotesize
	$$C= \displaystyle \max \,\left \lbrace      \left \lbrace   
	\displaystyle\frac{2}{1-\lambda}\,    \sin\left ( \displaystyle\frac{ \pi }{ N_b-1 }\right )\, \displaystyle \min_{0 \leq j \leq N_b-1 }  \left | \sin\left ( \displaystyle\frac{ \pi\,(2\,j+1 )}{  N_b-1 }\right )   \right| - \displaystyle  \frac{ 2\,\pi   }{  N_b \,(N_b-1)}  \,
	\displaystyle  \frac{ 1 }{   \lambda \, N_b  -1}  \right \rbrace  ,\eta_{ {\mathcal W} }        \right \rbrace $$
	\normalsize

	
	
	
	
	
	
	\noindent and consider the subdivision of the interval~$[0,1[$ into:
	
	$$N_m=\displaystyle \frac{1}{L_m}=(N_b-1)\, N_b^m$$
	
	\noindent sub-intervals of length~$L_m$. One has to determine a natural integer~$\tilde{N}_m$ such that the graph of~$\Gamma_{\mathcal W}$ on~$[0,1[$ can be covered by~$N_m \times \tilde{N}_m $ squares of side~\mbox{$ L_m $}. By considering, the vertical amplitude of the graph, one gets:

	$$\tilde{N}_m =\left \lfloor
	\displaystyle  \frac{C\, L_m^{2-D_{\mathcal W}   }}{L_m}\right\rfloor+1 
	\quad \text{i.e.} \quad \tilde{N}_m = \left \lfloor
	 C\,  L_m^{1-D_{\mathcal W}   }  \right\rfloor+1.
	$$
	
	\noindent  Thus:
	
	$$N_m \times \tilde{N}_m  =\displaystyle \frac{\tilde{N}_m}{L_m}=
	\displaystyle \frac{1}{L_m}\,\left \lfloor
	C\,  L_m^{1-D_{\mathcal W}   } \right \rfloor  + \displaystyle \frac{1}{L_m}.
	$$ 
	 
	\noindent The integer~$N_m \times \tilde{N}_m$ then obeys a power law of the form
	
	$$N_m \times \tilde{N}_m \approx c\, L_m^{-D_{\mathcal W}   }$$
	
	\noindent where~$c$ denotes a strictly positive constant.
	\end{proof}
\vskip 1cm

\noindent \textbf{Proof of Theorem 1.9.}\\

\noindent \emph{i}. \underline{Preliminary computations.}\\ 
	
 	For any pair of integers~$(i_m,j)$ of~\mbox{$\left \lbrace  0,\hdots,N_b-2  \right \rbrace^2$}:

	$$ \begin{array}{ccc} T_{i_m} \left (P_{j } \right)& =&\left( \displaystyle \frac{x_j+i_m}{N_b}, \lambda\, y_j + \cos\left(2\,\pi \,\left(\frac{x_j+i_m}{N_b}\right)\right) \right).
	\end{array}$$
	
  	 For any triplet of integers~$(i_m,i_{m-1},j)$ of~\mbox{$\left \lbrace  0,\hdots,N_b-2  \right \rbrace^3$}:

\footnotesize
$$ \begin{array}{ccc} &&T_{i_{m-1}}  \left ( T_{i_m}\left (P_{j } \right)\right)  =\qquad\qquad\qquad\qquad\qquad\qquad\qquad\qquad\qquad\qquad\qquad\qquad \qquad\qquad\qquad \\ \\
&&=\left( \displaystyle \frac{\frac{x_j+{i_m}}{N_b}+i_{m-1}}{N_b},
	\lambda^2\, y_j + \lambda\,\cos\left(2\,\pi \,\left(  \frac{x_j+{i_m}}{N_b} \right)\right)+ \cos\left(2\,\pi \,\left(\frac{\frac{x_j+{i_m}}{N_b}+i_{m-1} }{N_b}\right)\right) \right)  \qquad\\ \\
	&  &=\left( \displaystyle  \frac{x_j+{i_m}}{ N_b^2}+ \displaystyle \frac{  i_{m-1}}{N_b},
	\lambda^2\, y_j  +  \lambda \cos\left(2\,\pi \,\left(  \frac{x_j+{i_m}}{N_b} \right)\right)+  \cos\left(   2\,\pi \,\left( \displaystyle  \frac{x_j+{i_m}}{ N_b^2}+ \displaystyle \frac{  i_{m-1}}{N_b}\right)\right) \right).\\ \ 
	\end{array}$$
\normalsize
 	 For any quadruplet of integers~$(i_m,i_{m-1},i_{m-2},j)$ of~\mbox{$\left \lbrace  0,\hdots,N_b-2  \right \rbrace^4$}:

	\footnotesize
$$\begin{array}{ccc}  T_{i_{m-2}} \left (T_{i_{m-1}} \left ( T_{i_m}\left (P_{j } \right)\right)\right)  = \qquad\qquad\qquad\qquad\qquad\qquad\qquad\qquad\qquad\qquad\qquad\qquad\qquad\qquad \qquad\qquad \qquad&&\\  \\
=	\bigg(  \displaystyle  \frac{x_j+{i_m}}{ N_b^3}+ \displaystyle \frac{  i_{m-1}}{N_b^2} + \displaystyle  \frac{  i_{m-2}}{N_b},
\qquad \qquad\qquad\qquad\qquad\qquad\qquad\qquad\qquad\qquad\qquad\qquad\qquad\qquad\qquad \qquad\qquad\\
	   \lambda^3\, y_j \!+ \! \lambda^2\,  \cos\left(2\,\pi \,\left(  \frac{x_j+ i_m }{N_b} \right)\right)  
 \!+\!\lambda   \cos\left(2\,\pi \,\left( \displaystyle  \frac{x_j+{i_m}}{ N_b^2}+ \displaystyle \frac{  i_{m-1}}{N_b}\right)\right)
	\!+\!\cos\left(2\,\pi \,\left( \displaystyle  \frac{x_j+{i_m}}{ N_b^3}+ \displaystyle \frac{  i_{m-1}}{N_b^2} +  \displaystyle  \frac{  i_{m-2}}{N_b}\right)  \right)  \bigg). \!\!&&
	\end{array}$$
\normalsize

 	 Given a strictly positive integer~$m$, and two points~$X$ and~$Y$ of~$V_m$ such that:
	$$X  \underset{m }{\sim}  Y 
$$
	
 \noindent	 there exists a word~$\mathcal M$ of length~\mbox{$|{\mathcal M}|=m$}, on the graph~$\Gamma_{ \mathcal W} $, and an integer~$j$ of~\mbox{$\left \lbrace  0,\hdots,N_b-2  \right \rbrace^2$}, such that:
	
$$ X= T_{\mathcal M}  \left (P_{j } \right) \quad   , \quad Y= T_{\mathcal M}  \left (P_{j+1 } \right).
$$

 	 Let us write~$T_{\mathcal M}$ under the form:
	
$$ T_{\mathcal M}=T_{i_m} \circ T_{i_{m-1}} \circ \hdots \circ T_{i_1} 
$$

\noindent 	 where~\mbox{$(i_1,\hdots,i_m)\,\in\, \left \lbrace  0,\hdots,N_b-1  \right \rbrace^m$}.\\
	
	\newpage
 	 One has then:
	
$$ x \left (T_{\mathcal M}  \left (P_{j } \right) \right )  =
	\displaystyle  \frac{x_j }{ N_b^m}+ \displaystyle \sum_{k=1}^m\frac{  i_{k}}{N_b^k}  \quad , \quad
	x \left (T_{\mathcal M}  \left (P_{j+1 } \right) \right) =
	\displaystyle  \frac{x_{j+1} }{ N_b^m}+ \displaystyle \sum_{k=1}^m\frac{  i_{k}}{N_b^k}
$$

 \noindent	 and:

	$$ \left \lbrace \begin{array}{ccc}  y\left (T_{\mathcal M}  \left (P_{j } \right) \right) & =&\lambda^m\, y_{j } + \displaystyle \sum_{k=1}^{m } \lambda^{m-k}\,  \cos\left(2\,\pi \,\left(   \displaystyle  \frac{x_{j  } }{ N_b^{k}}+
	\displaystyle \sum_{\ell=0}^{k}\frac{  i_{ m-\ell}}{N_b^{k- \ell}}\right)\right)
	\\ \\
	y \left (T_{\mathcal M}  \left (P_{j+1 } \right) \right) & =&\lambda^m\, y_{j+1 } + \displaystyle \sum_{k=1}^{m } \lambda^{m-k}\,  \cos\left(2\,\pi \,\left(   \displaystyle  \frac{x_{j+1 } }{ N_b^{k}}+
	\displaystyle \sum_{\ell=0}^{k}\frac{  i_{ m-\ell}}{N_b^{k- \ell}}\right)\right)
	\end{array} \right.  .$$

	\normalsize

	\vskip 1cm
	
	\noindent \emph{ii}. \underline{Determination of the lower bound.}\\ 
	
\noindent Let us note that:
	
	\footnotesize
	$$\begin{array}{ccc}    h_{j,m} -\lambda^m\, \left (y_{j+1 }-y_{j } \right) =\qquad \qquad\qquad\qquad \qquad\qquad\qquad\qquad \qquad\qquad\qquad \qquad \qquad&   &\\ \\
=	\displaystyle \sum_{k=1}^{m } \lambda^{m-k}\, \left \lbrace  \cos\left(2\,\pi \,\left(   \displaystyle  \frac{x_{j+1 } }{ N_b^{k}}+
	\displaystyle \sum_{\ell=0}^{k}\frac{  i_{ m-\ell}}{N_b^{k- \ell}}\right)\right) -\cos\left(2\,\pi \,\left(   \displaystyle  \frac{x_{j  } }{ N_b^{k}}  -
	\displaystyle \sum_{\ell=0}^{k}\frac{  i_{ m-\ell}}{N_b^{k- \ell}}\right)\right) \right \rbrace \qquad &&\\
	\\  
	  =    -2 \displaystyle \sum_{k=1}^{m } \lambda^{m-k}\, \sin \left( \pi \,\left(   \displaystyle  \frac{x_{j+1 }-x_j }{  N_b^{k}} \right) \right) \, \sin \left(2\,\pi \,\left(   \displaystyle  \frac{x_{j+1 }+x_j }{ 2\,N_b^{k}}+
	\displaystyle \sum_{\ell=0}^{k}\frac{  i_{ m-\ell}}{N_b^{k- \ell}}\right) \right) . \qquad \qquad&&  \\
	\\ \\
	\end{array}$$
	\normalsize

	   Taking into account:

	\footnotesize
	$$\begin{array}{ccc}   \lambda^m\, \left (y_{j+1 }-y_{j } \right) & =  &
	\displaystyle\frac{\lambda^m}{1-\lambda}\, \left ( \cos\left ( \displaystyle\frac{2\,\pi\,(j+1)}{N_b-1}\right )- \cos\left ( \displaystyle\frac{2\,\pi\,j}{N_b-1}\right ) \right) \\
	   & =  &         -2\,  \displaystyle\frac{\lambda^m}{1-\lambda}\,   \sin\left ( \displaystyle\frac{2\,\pi\,(j+1-j)}{2\,(N_b-1)}\right )\, \sin\left ( \displaystyle\frac{2\,\pi\,(j+1+j)}{2\, (N_b-1)}\right ) \\
	& =  &         -2\,  \displaystyle\frac{\lambda^m}{1-\lambda}\,   \sin\left ( \displaystyle\frac{ \pi }{ N_b-1 }\right )\, \sin\left ( \displaystyle\frac{ \pi\,(2\,j+1 )}{  N_b-1 }\right ) \\
	\end{array}$$

\normalsize
	
	
	\footnotesize
	
	
	\normalsize

	\footnotesize
	
	\normalsize


	 
	
	\footnotesize


	\normalsize
	

	
	\normalsize

	

	
	\newpage
	
\noindent triangular inequality leads then to: 
	 \footnotesize

	 $$\begin{array}{ccc}   \left |  y \left (T_{\mathcal M}  \left (P_{j +1 } \right) \right)- y \left (T_{\mathcal M}  \left (P_{j   } \right)\right) \right|=\qquad\qquad\qquad\qquad\qquad\qquad\qquad\qquad \qquad\qquad\qquad \qquad\qquad\qquad
	 \qquad\qquad\qquad \qquad\qquad&&\\  \\
	 =\left | 
	   \lambda^m\, \left (y_{j +1 }-y_{j } \right)  -2\,\displaystyle \sum_{k=1}^{m } \lambda^{m-k}\, \sin \left(   \displaystyle  \frac{ \pi   }{  N_b^{k+1} \,(N_b-1)} \right)   \, \sin \left(   \displaystyle  \frac{ \pi \,(2\,j+1) }{  N_b^{k+1}\,(N_b-1)}+2\, \pi \,
	   \displaystyle \sum_{\ell=0}^{k}\frac{  i_{ m-\ell}}{N_b^{k- \ell}}\right)  
	   \right|\qquad\qquad\qquad \qquad\qquad \qquad&& \\  \\
	   \geq  \left | 
	   \lambda^m\,  \left | \left (y_{j +1 }-y_{j } \right)    \right| -2\,  \displaystyle \sum_{k=1}^{m } \lambda^{m-k}\, \left | \sin \left(   \displaystyle  \frac{ \pi   }{  N_b^{k+1} \,(N_b-1)} \right)   \, \sin \left(   \displaystyle  \frac{ \pi \,(2\,j+1) }{  N_b^{k+1}\,(N_b-1)}+2\, \pi \,
	   \displaystyle \sum_{\ell=0}^{k}\frac{  i_{ m-\ell}}{N_b^{k- \ell}}\right)  
	   \right|  \right| \qquad\qquad\qquad \qquad\qquad \quad&& \\ \\
	   =  \lambda^m\, \biggl | 
	     \displaystyle\frac{2}{1-\lambda}\,   \sin\left ( \displaystyle\frac{ \pi }{ N_b-1 }\right )\, \left | \sin\left ( \displaystyle\frac{ \pi\,(2\,j+1 )}{  N_b-1 }\right )   \right| \qquad\qquad\qquad \qquad\qquad\qquad\qquad \qquad\qquad\qquad\qquad
	     \qquad\qquad\qquad \qquad\qquad \qquad&&\\
	     -2\,  \displaystyle \sum_{k=1}^{m } \lambda^{ -k}\, \sin \left(   \displaystyle  \frac{ \pi   }{  N_b^{k+1} \,(N_b-1)} \right)   \, \left | \sin \left(   \displaystyle  \frac{ \pi \,(2\,j+1) }{  N_b^{k+1}\,(N_b-1)}+2\, \pi \,
	   \displaystyle \sum_{\ell=0}^{k}\frac{  i_{ m-\ell}}{N_b^{k- \ell}}\right)  
	   \right|  \biggl|.  \qquad\qquad\qquad\qquad \qquad\qquad\qquad\qquad \qquad\qquad\quad&&
	 \end{array}	$$
	 
	 \vskip 1cm
	 
	 \normalsize
	One also has:
	 \footnotesize

	 \footnotesize
$$\begin{array}{cccc}     2\,  \displaystyle \sum_{k=1}^{m } \lambda^{ -k}\, \sin \left(   \displaystyle  \frac{ \pi   }{  N_b^{k+1} \,(N_b-1)} \right)    \left | \sin \left(   \displaystyle  \frac{ \pi \,(2\,j+1) }{  N_b^{k+1}\,(N_b-1)}+2\, \pi \,
	 \displaystyle \sum_{\ell=0}^{k}\frac{  i_{ m-\ell}}{N_b^{k- \ell}}\right)  
	 \right| \leq    \qquad && \\  \\
	   \leq 2\,  \displaystyle \sum_{k=1}^{m } \lambda^{ -k}\, \sin \left(   \displaystyle  \frac{ \pi   }{  N_b^{k+1} \,(N_b-1)} \right)    \qquad   \qquad\qquad\qquad\qquad\qquad\qquad\qquad\qquad\qquad\quad && \\
	      \leq    2\,  \displaystyle \sum_{k=1}^{m } \lambda^{ -k}\,    \displaystyle  \frac{ \pi   }{  N_b^{k+1} \,(N_b-1)}    
	      \qquad\qquad\qquad\qquad\qquad\qquad \qquad\qquad\qquad\qquad\qquad\qquad&& \\
     =    \displaystyle  \frac{ 2\,\pi   }{  N_b \,(N_b-1)}  \,
  \displaystyle \sum_{k=1}^{m }     \displaystyle  \frac{ 1 }{  \lambda^{  k}\,N_b^{k }  } 
   \qquad\qquad\qquad\qquad\qquad\qquad \qquad\qquad\qquad \qquad\qquad\qquad&& \\
     =     \displaystyle  \frac{ 2\,\pi   }{  N_b \,(N_b-1)}  \,
  \lambda^{ -1}\, N_b^{-1 }     \displaystyle  \frac{ 1-\lambda^{ -m}\, N_b^{-m } }{  1-\lambda^{ -1}\, N_b^{-1 }  }   \qquad\qquad\qquad\qquad  \qquad\qquad \qquad\qquad\qquad&& \\ \\
  \leq     \displaystyle  \frac{ 2\,\pi   }{  N_b \,(N_b-1)}  \,
  \lambda^{ -1}\, N_b^{-1 }     \displaystyle  \frac{ 1 }{  1-\lambda^{ -1}\, N_b^{-1 }  }   \qquad\qquad\qquad\qquad\qquad \qquad\qquad\qquad \qquad&&  \\ \\
  =    \displaystyle  \frac{ 2\,\pi   }{  N_b \,(N_b-1)}  \,
   \displaystyle  \frac{ 1 }{   \lambda \, N_b  -1}  \qquad\qquad\qquad\qquad\qquad\qquad \qquad\qquad\qquad \qquad\qquad\qquad&&\\
\end{array}	$$
	 
	 \normalsize
	 
\noindent which yields:
	 
	 \footnotesize
	 $$\begin{array}{ccc}  -  2\,  \displaystyle \sum_{k=1}^{m } \lambda^{ -k}\, \sin \left(   \displaystyle  \frac{ \pi   }{  N_b^{k+1} \,(N_b-1)} \right)   \, \left | \sin \left(   \displaystyle  \frac{ \pi \,(2\,j+1) }{  N_b^{k+1}\,(N_b-1)}+2\, \pi \,
	 \displaystyle \sum_{\ell=0}^{k}\frac{  i_{ m-\ell}}{N_b^{k- \ell}}\right)  
	 \right|    \geq  - \displaystyle  \frac{ 2\,\pi   }{  N_b \,(N_b-1)}  \,
	 \displaystyle  \frac{ 1 }{   \lambda \, N_b  -1}  . \\
	 \end{array}	$$

	 \normalsize
	 
	 \newpage
	 \noindent Thus:

	 \footnotesize
	 $$\begin{array}{ccc}  \displaystyle\frac{2}{1-\lambda}\,   \sin\left ( \displaystyle\frac{ \pi }{ N_b-1 }\right )\, \left | \sin\left ( \displaystyle\frac{ \pi\,(2\,j+1 )}{  N_b-1 }\right )   \right| \qquad\qquad\qquad\qquad \qquad\qquad \qquad
	 &&\\-  2\,  \displaystyle \sum_{k=1}^{m } \lambda^{ -k}\, \sin \left(   \displaystyle  \frac{ \pi   }{  N_b^{k+1} \,(N_b-1)} \right)   \, \left | \sin \left(   \displaystyle  \frac{ \pi \,(2\,j+1) }{  N_b^{k+1}\,(N_b-1)}+2\, \pi \,
	 \displaystyle \sum_{\ell=0}^{k}\frac{  i_{ m-\ell}}{N_b^{k- \ell}}\right)  
	 \right|    \geq  &&\\ \\
	  \geq \displaystyle\frac{2}{1-\lambda}\,   \sin\left ( \displaystyle\frac{ \pi }{ N_b-1 }\right )\, \left | \sin\left ( \displaystyle\frac{ \pi\,(2\,j+1 )}{  N_b-1 }\right )   \right| - \displaystyle  \frac{ 2\,\pi   }{  N_b \,(N_b-1)}  \,
	 \displaystyle  \frac{ 1 }{   \lambda \, N_b  -1}  . \qquad\qquad&& \\
	 \end{array}	$$

	 \normalsize
	 \vskip 1cm
	 
	 \begin{lemma}The following results hold:

	 $$\sin\left ( \displaystyle\frac{ \pi\,(2\,j+1 )}{  N_b-1 }\right ) =0 \quad \text{if and only if~$N_b$ is even and~$j=\displaystyle\frac{  N_b}{  2 } -1 $.}$$

	 \end{lemma}
 
 \vskip 1cm

\noindent \textbf{Proof of lemma 2.3.}\\

	 \noindent Since~$0 \leq j \leq  N_b-1$, one has:

	 $$0 \leq 2\,  j +1\leq  2\, N_b-1\quad \text{and thus} \quad  0 \leq  \displaystyle\frac{  2\,j+1  }{  N_b-1 } \leq  2 +\displaystyle\frac{  1 }{  N_b-1 } .$$
	 
	 \noindent Then,~$\sin\left ( \displaystyle\frac{ \pi\,(2\,j+1 )}{  N_b-1 }\right ) =0$ in the sole cases:
	 
	 $$  \displaystyle\frac{  2\,j+1  }{  N_b-1 } = 1 \quad \text{or} \quad  \displaystyle\frac{  2\,j+1  }{  N_b-1 } = 2.
	 	 $$
	 	 
	 	 \noindent The second case has to be rejected, since it would lead to:
	 	 
	 	  $$  j=N_b-\displaystyle\frac{  3}{  2 }   \, \notin \, \N
	 	 $$
	 	 
	 	 \noindent The only possibility is thus when~$N_b$ is an even number:
	 	 $$  j=\displaystyle\frac{  N_b}{  2 } -1  $$
	 \noindent The converse side is obvious. \\
	 
	 \noindent \hfill $\square$
	 
	 \vskip 1cm
	 \noindent \underline{First case:}~$\sin\left ( \displaystyle\frac{ \pi\,(2\,j+1 )}{  N_b-1 }\right ) \neq 0$.
	 \noindent  One has then:
	 
	 $$\left | \sin\left ( \displaystyle\frac{ \pi\,(2\,j+1 )}{  N_b-1 }\right )   \right|  \geq 
	 \displaystyle \min_{0 \leq j \leq N_b-1}\, \left | \sin\left ( \displaystyle\frac{ \pi\,(2\,j+1 )}{  N_b-1 }\right )   \right| 
	 =\left | \sin\left ( \displaystyle\frac{ \pi}{  N_b-1 }\right )   \right|\geq  \displaystyle\frac{ 2}{  N_b-1 }.$$
	 
	 \noindent This leads to:
	 
	  \footnotesize
	 $$\begin{array}{ccc}  \displaystyle\frac{2}{1-\lambda}\,   \sin\left ( \displaystyle\frac{ \pi }{ N_b-1 }\right )\, \left | \sin\left ( \displaystyle\frac{ \pi\,(2\,j+1 )}{  N_b-1 }\right )   \right| \qquad\qquad\qquad\qquad \qquad\qquad \qquad \qquad&&\\
	 -  2\,  \displaystyle \sum_{k=1}^{m } \lambda^{ -k}\, \sin \left(   \displaystyle  \frac{ \pi   }{  N_b^{k+1} \,(N_b-1)} \right)   \, \left | \sin \left(   \displaystyle  \frac{ \pi \,(2\,j+1) }{  N_b^{k+1}\,(N_b-1)}+2\, \pi \,
	 \displaystyle \sum_{\ell=0}^{k}\frac{  i_{ m-\ell}}{N_b^{k- \ell}}\right)  
	 \right|    \geq  \qquad &&\\ \\
	 \geq \displaystyle\frac{2}{1-\lambda}\,   \sin\left ( \displaystyle\frac{ \pi }{ N_b-1 }\right )\,   \sin\left ( \displaystyle\frac{ \pi   }{  N_b-1 }\right )   - \displaystyle  \frac{ 2\,\pi   }{  N_b \,(N_b-1)}  \,
	 \displaystyle  \frac{ 1 }{   \lambda \, N_b  -1}   \qquad\qquad\qquad \qquad&& \\ \\
	 \geq \displaystyle\frac{2}{1-\lambda}\,    \displaystyle\frac{ 4}{ (N_b-1 )^2 }       - \displaystyle  \frac{ 2\,\pi   }{  N_b \,(N_b-1)}  \,
	 \displaystyle  \frac{ 1 }{   \lambda \, N_b  -1}  \qquad\qquad\qquad\qquad\qquad  \qquad\qquad \qquad&& \\ \\
	  =   \displaystyle\frac{ 2}{ N_b-1 } \,   \left \lbrace \displaystyle\frac{4}{1-\lambda}\,     \displaystyle\frac{ 1}{  N_b-1 }       - \displaystyle  \frac{  \pi   }{  N_b  }  \,
	 \displaystyle  \frac{ 1 }{   \lambda \, N_b  -1}  \right \rbrace  \qquad\qquad\qquad\qquad \qquad\qquad\qquad \qquad&& \\ \\
	  =   \displaystyle\frac{ 2}{ N_b\,(N_b-1 )\,(1-\lambda)\,( \lambda\,N_b-1)} \,  
	   \left \lbrace 4\,N_b\, ( \lambda \, N_b  -1)      -   \pi  \,(1-\lambda)\,  (N_b-1 )\right \rbrace . \qquad\qquad && \\
	 \end{array}	$$

	 \normalsize

	\noindent  The function
	 
	 $$\lambda \mapsto 4\,N_b\, ( \lambda \, N_b  -1)      -   \pi  \,(1-\lambda)\,  (N_b-1 )$$
	 
	 \noindent is affine and strictly increasing in~$\lambda$, and quadratic and strictly increasing in~$N_b$, for strictly positive values of~$N_b$. This ensures the positivity of:
	 
	 \footnotesize
	 $$\begin{array}{ccc}  \displaystyle\frac{2}{1-\lambda}\,   \sin\left ( \displaystyle\frac{ \pi }{ N_b-1 }\right )\, \left | \sin\left ( \displaystyle\frac{ \pi\,(2\,j+1 )}{  N_b-1 }\right )   \right| \qquad\qquad\qquad\qquad\qquad\qquad\qquad  && \\
	 -  2\,  \displaystyle \sum_{k=1}^{m } \lambda^{ -k}\, \sin \left(   \displaystyle  \frac{ \pi   }{  N_b^{k+1} \,(N_b-1)} \right)   \, \left | \sin \left(   \displaystyle  \frac{ \pi \,(2\,j+1) }{  N_b^{k+1}\,(N_b-1)}+2\, \pi \,
	 \displaystyle \sum_{\ell=0}^{k}\frac{  i_{ m-\ell}}{N_b^{k- \ell}}\right)  
	 \right|     . &&\\
	 \end{array}	$$
	  
	   \normalsize
	   
	   \vskip 1cm

\noindent \underline{Second case:}~$\sin\left ( \displaystyle\frac{ \pi\,(2\,j+1 )}{  N_b-1 }\right ) = 0$.\\

\noindent One has then:

	   \footnotesize
	   
	   $$\begin{array}{ccc}   \left |  y \left (T_{\mathcal M}  \left (P_{j +1 } \right) \right)- y \left (T_{\mathcal M}  \left (P_{j   } \right)\right) \right| \geq  2\, \lambda^m\, \left |    \displaystyle \sum_{k=1}^{m } \lambda^{ -k}\, \sin \left(   \displaystyle  \frac{ \pi   }{  N_b^{k+1}  } \right)   \, \left | \sin \left(   \displaystyle  \frac{ \pi  }{  N_b^{k+1}}+2\, \pi \,
	   \displaystyle \sum_{\ell=0}^{k}\frac{  i_{ m-\ell}}{N_b^{k- \ell}}\right)  
	   \right|  \right| && \\ 
	   \end{array}	$$

	   \vskip 1cm
	   
	   \normalsize
	   
	   \noindent Thanks to the periodic properties of the sine function, one may only consider the case when:

	   $$0 \leq \displaystyle  \frac{ \pi }{  N_b^{k+1} }+2\, \pi \,
	   \displaystyle \sum_{\ell=0}^{k}\frac{  i_{ m-\ell}}{N_b^{k- \ell}} \leq \displaystyle  \frac{ \pi }{2}.$$

	   \noindent Thus:
	   
	   $$\begin{array}{ccc}   \left |  y \left (T_{\mathcal M}  \left (P_{j +1 } \right) \right)- y \left (T_{\mathcal M}  \left (P_{j   } \right)\right) \right| &\geq & 
	   \displaystyle \sum_{k=1}^{m } \lambda^{ -k}\,   \displaystyle  \frac{2   }{  N_b^{k+1}  }    \, 
	   \left \lbrace   \displaystyle  \frac{2  }{  N_b^{k+1} }+2\, 
	   \displaystyle \sum_{\ell=0}^{k}\frac{  i_{ m-\ell}}{N_b^{k- \ell}}\right \rbrace  \\
	     & 
	      \geq & \displaystyle \sum_{k=1}^{m } \lambda^{ -k}\,  \displaystyle  \frac{2   }{  N_b^{k+1}  }    \, 
	     \left \lbrace   \displaystyle  \frac{2  }{  N_b^{k+1} } \right \rbrace  \\
	     &=&  \displaystyle  \frac{4 \,\lambda^{ -1}   }{  N_b^{4}  }    \,   \displaystyle  \frac{1- \lambda^{ -m}\,N_b^{ -2m}  }{  1-\lambda^{ -1}\,N_b^{ -2} } \\ 
	     & =&\displaystyle  \frac{4   }{  N_b^{2} \,(N_b-1)^2}    \,   \displaystyle  \frac{1- \lambda^{ -m}\,N_b^{ -m}  }
	     {   \lambda \,N_b-1}  \\
	     &=&   \displaystyle  \frac{4   }{  N_b^{2}  }    \,   \displaystyle  \frac{1-  N_b^{ -2}  }
	     {    N_b^2-1  }  .
	      \\
	   \end{array}	$$

	   \normalsize

  \vskip 1cm
  
  \noindent \underline{General case:}~the above results enable us to obtain the predominant term of the lower bound of~
  \mbox{$\left |  y \left (T_{\mathcal M}  \left (P_{j +1 } \right) \right)- y \left (T_{\mathcal M}  \left (P_{j   } \right)\right) \right|$}, which is thus:
	
$$ \lambda^{m  }
	=  e^{ m\,(D_{\mathcal W}-2)\, \ln N_b}= N_b^{ m\,(D_{\mathcal W}-2)}  = L_m^{2-D_{\mathcal W}}  \,(N_b-1)^{2-D_{\mathcal W}}  . $$

	 \noindent \emph{ii.} \underline{Determination of the upper bound}.\\
	
	 \noindent One has:

	\footnotesize
	
$$\begin{array}{ccc}  \left | h_{j,m}\right| & \leq  &
	\displaystyle\frac{2\,\lambda^m }{1-\lambda}\,  \displaystyle\frac{ \pi^2\,(2\,j+1)}{(N_b-1)^2}   +
	2\,\displaystyle \sum_{k=1}^{m } \lambda^{m-k}\, \pi \,\left\lbrace   \displaystyle \frac{2\,j+1}{(N_b-1)\, N_b^k}+
	2\, \displaystyle \sum_{\ell=0}^{k}\frac{  i_{ m-\ell}}{N_b^{k- \ell}}\right \rbrace \,   \displaystyle  \frac{  \pi }{ (N_b-1)\, N_b^k}
	\\ \\
	& =&       \displaystyle\frac{2\,\lambda^m}{1-\lambda}\,  \displaystyle\frac{ \pi^2\,(2\,j+1)}{(N_b-1)^2}  +
	\displaystyle  \frac{  2\,\pi^2\, \lambda^m  }{  N_b-1 }\, \displaystyle \sum_{k=1}^{m }  \,\left\lbrace   \displaystyle \frac{(2\,j+1)\, \lambda^{ -k}}{(N_b-1)\, N_b^{2k}}+
	2\, \displaystyle \sum_{\ell=0}^{k}\frac{  i_{ m-\ell}\,  \lambda^{ -k} }{N_b^{2k- \ell}}\right \rbrace
	\\ \\
	& = &      \displaystyle\frac{2\,\lambda^m}{1-\lambda}\,  \displaystyle\frac{ \pi^2\,(2\,j+1)}{(N_b-1)^2} \\
	&& +
	\displaystyle  \frac{  2\,\pi^2\, \lambda^m }{  N_b-1 }\,   \left\lbrace
	\displaystyle \frac{\lambda^{-1}\,N_b^{-2}\,(2\,j+1) }{(N_b-1) } \,  \displaystyle \frac{  (1- \lambda^{ -m }\,N_b^{ -2m })}{1- \lambda^{ - 1}\,N_b^{ - 2}} +
	2\, \displaystyle \sum_{k=1}^m \frac{  (N_b-1)\,  \lambda^{ -k}
	}{N_b^{2k }} \displaystyle \frac{1-N_b^{-k-1}}{1-N_b^{-1}}\right \rbrace \\ \\
	& \leq  &      \displaystyle\frac{2\,\lambda^m}{1-\lambda}\,  \displaystyle\frac{ \pi^2\,(2\,N_b-1)}{(N_b-1)^2}
	+
	\displaystyle  \frac{  2\,\pi^2\, \lambda^m }{  N_b-1 }\,
	\displaystyle \frac{ (2\,N_b-1) }{(N_b-1) } \,  \displaystyle \frac{  (1- \lambda^{ -m }\,N_b^{ -2m })}{  \lambda \,N_b^{   2}-1}   \\ \\
	&& +
	\displaystyle  \frac{  2\,\pi^2\, \lambda^m }{  N_b-1 }\,
	2\, \displaystyle  \frac{ \lambda^{ -1 }\,N_b^{ -2  }\, (N_b-1)\, (1- \lambda^{ -m }\,N_b^{ -2m })
	}{ (1-N_b^{-1})\,(1-\lambda^{ -1 }\,N_b^{ -2  })}  \\ \\
	&& -
	\displaystyle  \frac{  2\,\pi^2\, \lambda^m }{  N_b-1 }\,
	2\, \displaystyle  \frac{ \lambda^{ -1 }\,N_b^{ -3  }\, (N_b-1)\, (1- \lambda^{ -m }\,N_b^{ -3m })}{ (1-N_b^{-1})\,(1-\lambda^{ -1 }\,N_b^{ -3  })}
	\\ \\
	& \leq  &      \displaystyle\frac{2\,\lambda^m}{1-\lambda}\,  \displaystyle\frac{ \pi^2\,(2\,N_b-1)}{(N_b-1)^2}  +
	\displaystyle  \frac{  2\,\pi^2\, \lambda^m }{  N_b-1 }\,
	\displaystyle \frac{ (2\,N_b-1) }{(N_b-1) } \,  \displaystyle \frac{ 1}{  \lambda \,N_b^{   2}-1}   \\ \\
	&& +
	\displaystyle  \frac{  4\,\pi^2\, N_b\, \lambda^m }{   N_b-1  }\,\left \lbrace   \displaystyle  \frac{ 1}{   \lambda \,N_b^{  2  }-1  }  -
	\displaystyle  \frac{ 1}{   \lambda \,N_b^{  3  }-1  }\right \rbrace
	\\ \\
	& =  &2\, \pi^2\,\lambda^m \,\left \lbrace
	\displaystyle  \frac{   (2\,N_b-1)\, \lambda\, (N_b^2-1) } {(N_b-1)^2 \, (1- \lambda )\,(\lambda \,N_b^{   2}-1) }    +
	\displaystyle  \frac{  2\, N_b } {   (\lambda \,N_b^{ 2  }-1)\, (\lambda \,N_b^{ 3  }-1)  } \right \rbrace .
	\\ \\
	
	\end{array}$$
	 
 \normalsize
	
	 Since:

$$x\left ( T_{\mathcal M}  \left (P_{j+1 } \right) \right)-x\left ( T_{\mathcal M}  \left (P_{j } \right)\right)=
	\displaystyle \frac{1}{(N_b-1)\, N_b^m} $$
	
 \noindent	 and:

$$ D_{\mathcal W}= 2+\displaystyle \frac{\ln \lambda}{\ln N_b} \quad , \quad \lambda= e^{(D_{\mathcal W}-2)\, \ln N_b}= N_b^{(D_{\mathcal W}-2) } $$
	
	\newpage
 \noindent	 one has thus:

	\footnotesize

$$ \begin{array}{ccc}  \left | h_{j,m}\right| & \leq  & 2\, \pi^2\,L_m^{2-D_{\mathcal W}   }\, \left (N_b-1\right)^{2-D_{\mathcal W}   }\,\,\left \lbrace
	\displaystyle  \frac{   (2\,N_b-1)\, \lambda\, (N_b^2-1) } {(N_b-1)^2 \, (1- \lambda )\,(\lambda \,N_b^{   2}-1) }    +
	\displaystyle  \frac{  2\, N_b } {   (\lambda \,N_b^{ 2  }-1)\, (\lambda \,N_b^{ 3  }-1)  } \right \rbrace .
	\\ \\
	\end{array}$$
\normalsize

	\noindent \hfill $\square$

\vskip 1cm

\begin{remark}
	
	 In~\cite{Hunt1998}, B.~Hunt uses the fact that the Hausdorff dimension of a fractal set~$\mathcal F$ can be obtained by means of what is called the~$t-$energy,~\mbox{$t\,\in\,\R$}, of a Borel measure supported on~$\mathcal F$ (one may refer to~\cite{Falconer1985}, for instance):

$$ I_t(\mu)= \displaystyle \int \!\!\!\int \displaystyle \frac{d\mu(x)\,d\mu(x')}{|x-x'|^t}$$
	
 \noindent	 which enables one to obtain:

$$ \dim {\mathcal F}=\displaystyle \sup\,\left \lbrace t \,\in\,\R, \, \text{$\mu$ supported on $\mathcal F$}  \big |\,  I_t(\mu)< + \infty \right \rbrace $$
	
	 A lower bound~$t_0$ of the Hausdorff dimension can thus be obtained by building a measure~$\mu$ supported on~$\mathcal F$ such that:

$$ I_{t_0}(\mu)< + \infty  .$$

	 B.~Hunt proceeds as follows: he introduces the measure~$\mu_{\mathcal W}$ supported on~$\Gamma_{\mathcal W}$, induced by the Lebesgue measure~$\mu$ on~$[0,1]$. Thus:

$$ I_t\left (\mu_{\mathcal W}\right )=
	\displaystyle \int \!\!\!\int \displaystyle \frac{d\mu_{\mathcal W}(x)\,d\mu_{\mathcal W}(x')}{\left \lbrace |x-x'|^2+ |{\mathcal W}(x)-{\mathcal W}(x')|^2\right \rbrace ^{\frac{t}{2}}}  .$$
	
	 We could also have used a similar argument since, in our case:

$$  |x- x' |^{2-D_{\mathcal W}}\lesssim  |{\mathcal W}(x)-{\mathcal W}(x')| \lesssim |x- x' |^{2-D_{\mathcal W}}.
$$

\end{remark}

\vskip 1cm

\centerline{\textbf{Thanks}}

\vskip 1cm
The author would like to thank the anonymous referee for his careful reading, and his very pertinent suggestions and advices, which helped a lot improving the original work.

\vskip 2cm
















\newpage
\bibliographystyle{alpha}
 
\bibliography{DimensionGammaW}

\end{document}